\documentclass[12pt,leqno]{article}
\newcommand{\EE}{{\rm I\kern-2pt E}}
\newcommand{\RR}{{\rm I\kern-2pt R}}
\newcommand{\DD}{{\rm I\kern-2pt D}}
\newcommand{\PP}{{\rm I\kern-2pt P}}
\newcommand{\NN}{{\rm I\kern-2pt N}}
\newcommand{\dd}{{\rm \kern 3pt I\kern-9pt d}}

\title{\large DIRICHLET FORMS METHODS : AN APPLICATION TO THE PROPAGATION OF THE ERROR DUE TO THE EULER SCHEME}
\author{\sc Nicolas Bouleau\\
{\tt bouleau@enpc.fr}}
\date{\it Ecole des Ponts, ParisTech}
\begin{document}
\maketitle

\noindent{\bf Abstract.} We present recent advances on Dirichlet forms methods either to extend financial models beyond the usual stochastic
calculus or to study stochastic models with less classical tools. In this spirit, we interpret the asymptotic error on the solution of an sde due 
to the Euler scheme (Kurtz and Protter [Ku-Pr-91a]) in terms of a Dirichlet form on the Wiener space, what allows to propagate this error thanks to
 functional calculus. \\


\noindent{\sf\large Introduction}\\

Considering a Dirichlet form amounts to consider a strongly continuous symmetric contraction  semi-group on an $L^2$-space which possesses
in addition the property of being positive on positive functions 
(cf. [Fu-Os-Ta-94], [Bo-Hi-91 ], [Ma-R\"o-92]). It is a particular case of Markovian 
potential theory, with several special features due to the use of Hilbertian techniques and to the fact that positivity and contraction properties
extend  to infinite dimensional framework thanks to Fatou's lemma in measure theory. Many Dirichlet structures are constructively obtained
on the Wiener space and on the fundamental spaces of probability theory (Poisson space, Monte Carlo space) which may be thought as
hypotheses in order to study error propagation through stochastic models (cf. [Bou-03b]).

Since the discovery by M. Fukushima, at the end of the seventies, that Dirichlet forms allow to extend the stochastic calculus to processes
which are not semi-martingales (cf. [Fuk-80]) a lot of works have been developed in this direction, even beyond the 
Dirichlet forms framework. To this extend we quote the approach to time-dependent Dirichlet forms developed by Oshima [Osh-92] and 
the more recent approach of Stannat [Sta-99] and Trutnau [Tru-00] about a new theory of generalized Dirichlet forms. As in finance the heart of the complete market property and more generally of the portfolio management is the 
stochastic integral, a particular interest has been devoted to methods giving rise to new stochastic integrals.

We shall give, at first, a short outlook on recent results related to Dirichlet forms and connected with financial motivations. We include some
Malliavin calculus approaches when they amount to the use of the Ornstein-Uhlenbeck structure on the Wiener space. After recalling, in a second
part, the main properties of Dirichlet forms and the interpretation of the functional calculus on the squared field operator 
in terms of error propagation, we focuse, in a third part, on the question of the asymptotic error due to the resolution of a stochastic differential
equation by the Euler scheme. We show that the asymptotic error may be represented by a Dirichlet structure on the Wiener space and we
apply this to propagate the error on the example of a level volatility model for pricing and hedging procedures. 
We put the general question of the validity of such a propagation as an {\it asymptotic calculus principle}, 
and we give partial arguments for this principle.\\

\noindent {\sf\large I. Some recent works.}\\

First must be mentionned the idea of using Malliavin's integration by parts technique to speed up the computation of the Greeks or other quantities
in finance. After the collective papers of Fournier and {\it al} [Fo-La-Le-Li-To-99] [Fo-La-Le-Li-01], improvements have been 
brought to complex options [Go-Ko-01] and to the more general question of
 the sensitivity to some parameters with the aim of calibration of a model. As integration by parts formulae are available in more general Dirichlet forms situations than the Ornstein-Uhlenbeck
structure on the Wiener space (cf. [Bou-03b] Chapter V), the same approach may be performed for instance on the Poisson space for
studying models with jumps [El-Pr-04].

One of the first success of Malliavin calculus was about proving existence of densities for solutions of sde's with smooth coefficients
and Dirichlet forms methods have been able to extend such results to the case of Lipschitz coefficients [Bo-Hi-91]. 
Several authors remarked that these tools give also means of improving the computation of densities and  establishing estimates for the laws
 of random variables with some regularity assumption. Let us quote ([Ko-Pe-02] , [Ca-Fe-Nu-98], 
[Bo-Ek-To-04], [Bou-05b]) whose results aren't limited to applications in finance. With suitable hypotheses it is possible,  to get explicit closed formulae for the density even with some liberty in the choice of a weight function
allowing an optimization for Monte Carlo simulation.

After the classical works of M. Fukushima and Y. Le Jan on stochastic calculus for additive functionals of symmetric 
Markov processes associated with a Dirichlet form [LeJ-78] the role of past and future
$\sigma$-algebras have been clarified by Lyons and Zheng (cf [Ly-Zh-98] [Tru-00]) and the main current of research, in order to leave
the semi-martingale context, starts with
the abstract definition of a Dirichlet process as sum of a local martingale and a process with zero quadratic variation (see [F\"ol-80]).  
 Because the quadratic variation, as formal Dirichlet form, does not possess the closedness property, the Dirichlet form
framework is replaced here by functional analytic arguments. The integral is  generally defined by a discretization procedure
(cf. [F\"{o}l-81], [Bou-85], [F\"o-Pr-Sh-95] ) or by a regularization procedure (see [Ru-Va-95], [Ru-Va-96]). These ways have been 
deepened with the center example of the fractional Brownian motion
(cf.  [Er-Ru-98],  [Zah-98], [Fe-LaP-99], [Al-Ma-Nu-00], [Ru-Va-00], [Gr-No-03], [Gr-Ru-Va-03]). The connections of these works with finance are many : attempting to generalize Girsanov theorem in order to define martingale
measures by erasing more general drifts and using generalized stochastic integration (forward, symmetric and backward integrals) in order
to deal with exotic models (cf [Fl-Ru-Wo-03]). About ``inside trading" and the use of forward integral it is worth to 
quote [Le-Na-Nu-03].

At last, let us mention some uses of Dirichlet forms or Malliavin calculus to deal with processes with jumps by equipping the 
general Poisson 
space with a differential structure (cf. [De-Gr-Po-99],  [Me-Pr-03],  ) and the forthcoming book of P. Malliavin and A. Thalmaier [Ma-Th-05] whose last chapter
is devoted to calculus of variations for markets with jumps, the other ones being strongly related with the above topics.\\

\noindent {\sf\large II. Dirichlet forms theory seen as error propagation theory.}\\

Let us begin with a very simple but crucial remark about the magnitude of errors. If we consider an erroneous quantity with a 
centered small error and apply to it a non linear map, we observe by an easy Taylor expansion argument that

- the error is no more centered in general : a bias appears

- the variance transmit with a first order calculus.

\noindent Now if we go on, applying anew several non-linear applications

- the variances and the biases keep (except special cases) the same order of magnitude

- the biases follow a second order differential calculus involving the variances.

\noindent With natural notation
$$\sigma^2_{n+1}=f^{\prime 2}_{n+1}(x_n)\sigma^2_n$$
$${\mbox{bias}_{n+1}}=f^\prime_{n+1}(x_n){\mbox{bias}_{n}}+\frac{1}{2}f^{\prime\prime}_{n+1}(x_n)\sigma^2_n.$$
The first relation has been discovered, even in several dimension with correlation between the errors, 
by Gauss at the beginning of the nineteenth century. 

From this observation, in order to represent the propagation of small errors we may consider that

1) the variances of errors have to be managed by a quadratic first order differential operator $\Gamma$,

2) the biases of errors have to be represented by a linear second order differential operator $A$,

\noindent the propagation of errors being the result of the following change of variable formulae :
$$\Gamma[F(X_1,\ldots,X_m),G(Y_1,\ldots,Y_n)]=\sum_{ij}F^\prime_i(X_1,\ldots,X_m)G^\prime_j(Y_1,\ldots,Y_n)\Gamma[X_i,Y_j]$$
$$A[F(X_1,\ldots,X_m)]=\sum_iF^\prime_i(X_1,\ldots,X_m)A[X_i]+\frac{1}{2}\sum_{ij}F^{\prime\prime}_{ij}(X_1,\ldots,X_m)\Gamma[X_i,X_j].$$
Because of these propagation rules for the variances and the biases, little errors may be thought as {\it second order} vectors. This old notion
of differential geometry has been revived at the beginning of the eighties by the study of semi-martingales on manifolds 
(cf. [Sch-82] [Mey-82] [Eme-89]).

Now, instead of germs of semi-martingales and second order vectors, we will use Dirichlet forms, carr\'e du champ and generator.
 There are two important
reasons for this, that I shall give just after recalling some  definitions and examples.\\

\noindent{\bf Definition} {\it An error structure is a term}\index{error structure}
$$
S=( \Omega ,\mathcal{A}, \PP ,\DD,\Gamma )
$$
{\it where $(\Omega ,\mathcal{A},\PP )$ is a probability space, and:}

\begin{enumerate}
\item[(1)] $\DD$ {\it is a dense subvector space of $L^2 (\Omega ,
\mathcal{A}, \PP )$ (also denoted $L^2 (\PP )$)};\index{$\DD$ domain of $\Gamma$}

\item[(2)] {\it $\Gamma$ is a positive symmetric bilinear application
from $\DD \times \DD$ into $L^1 (\PP )$ satisfying
``the functional calculus of class $\mathcal{C}^1 \cap \mathrm{Lip}$''.
This expression means}\index{functional calculus of class $\mathcal{C}^1 \cap \mathrm{Lip}$}\index{$\Gamma$
quadratic error operator}
$$
\forall u\in \DD^m , \quad \forall v\in \DD^n ,
\quad \forall F\colon \RR^m \to \RR ,\quad \forall G\colon
\RR^n \to \RR
$$
{\it with $F$, $G$ being of class $\mathcal{C}^1$ and Lipschitzian, we have $F(u) \in \DD$,
$G(v) \in \DD$ and}
$$
\Gamma [F(u),G(v)] =\sum_{i,j} \frac{\partial F}{\partial x_i} (u)
\frac{\partial G}{\partial x_j} (v)\Gamma \bigl[ u_i ,v_j\bigr] 
\quad \PP \hbox{-a.s.} ;
$$

\item[(3)] {\it the bilinear form $\mathcal{E} [u,v] =\frac{1}{2}
\EE \bigl[ \Gamma [u,v]\bigr]$ is ``closed"\index{closed form}. This means that
the space $\DD$ equipped with the norm}\index{$\mathcal{E}$ Dirichlet form}\index{Dirichlet form}
$$
\| u\|_{\DD} =\left( \| u\|^2_{L^2 (\PP )} +
\mathcal{E} [u,u] \right)^{1/2}
$$
{\it is complete}.
\end{enumerate}

{\it If, in addition}
\begin{enumerate}
\item[(4)] {\it the constant function $1$ belongs to $\DD$ (which implies
$\Gamma [1] =0$ by property 2), we say that the error structure is Markovian.}\index{error structure, Markovian}
\end{enumerate}

\noindent
We will always write $\mathcal{E} [u]$ for $\mathcal{E} [u,u]$ and
$\Gamma [u]$ for $\Gamma [u,u]$.

With this definition, the form $\mathcal{E}$ is known in the
literature as a {\it local Dirichlet form}\index{Dirichlet form} on $L^2 (\Omega ,\mathcal{A},
\PP )$ that possesses a ``squared field'' operator\index{squared field operator} (or a ``carr\'e du
champ'' operator)\index{carr\'e du champ operator} $\Gamma$. These notions are usually studied on
$\sigma$-finite measurable spaces. We limit ourselves herein to probability
spaces both for the sake of simplicity and because we will use images and products of
error structures.

Under very weak additional assumptions, to an error structure\index{error structure} (also to a Dirichlet form\index{Dirichlet
form} on a $\sigma$-finite measurable space)  a
strongly-continuous contraction semigroup\index{semigroup, strongly continuous} $\bigl( P_t\bigr)_{t\geq 0}$
on $L^2 (\PP)$ can be uniquely associated, which is symmetric with respect to $\PP$ and 
sub-Markov. This semigroup has a generator\index{generator} $(A,\mathcal{D} A)$, a
self-adjoint operator that satisfies:
$$
A\bigl[ F(u)\bigr] =\sum_i \frac{\partial F}{\partial x_i} (u)
A\bigl[ u_i\bigr] +\frac{1}{2} \sum_{i,j} \frac{\partial^2 F}{\partial x_i
\partial x_j} (u)\Gamma \bigl[ u_i,u_j\bigr] \ \ \PP \hbox{-a.s.}
$$
for $F\colon \RR^m \to \RR$ of class $\mathcal{C}^2$ with bounded
derivatives and $u\in (\mathcal{D} A)^m$ such that $\Gamma \bigl[ u_i
\bigr] \in L^2 (\PP )$.

\noindent
{\bf Example 1.}(Ornstein-Uhlenbeck structure in dimension 1)

$\Omega  = \RR$, 
$\quad\mathcal{A}  =  \hbox{Borel $\sigma$-field} \ \mathcal{B} (
\RR )$,
$\quad\PP = \mathcal{N}(0,1) \ \hbox{reduced normal law} $,
$\quad\DD  = H^1 \bigl( \mathcal{N}(0,1)\bigr)  = \bigl\{ u\in L^2
({\PP} ), u' $ in the distribution sense belongs to $ L^2 (\PP )\bigr\} $,
$\quad\Gamma [u]  =  u^{\prime 2}$,
then $\bigl( \RR ,\mathcal{B} (\RR ),
\mathcal{N}(0,1), H^1 (\mathcal{N}(0,1)), \Gamma \bigr)$ is an error structure with
 generator\index{generator}
$$
\mathcal{D} A = \bigl\{ f\in L^2 (\PP)\colon \hbox{$f'' -xf'$ in the distribution sense} \in L^2 (\PP)\bigr\}
$$
$$
Af =\frac{1}{2} f'' -\frac{1}{2} I\cdot f'
$$ 
where $I$ is the identity map on $\RR$.

\medskip
\noindent
{\bf Example 2.} (Monte Carlo structure in dimension 1)

$\Omega  = [0,1]$,
$\quad\mathcal{A} = \hbox{Borel $\sigma$-field}$,
$\quad\PP =\hbox{Lebesgue measure}$,
$\quad\DD  = \bigl\{ u\in L^2 \bigl( [0,1] ,dx\bigr)$
the derivative $u'$ in the distribution sense over $]0,1[$ belongs to 
$L^2 ([0,1], dx)\bigr\}$,
$\quad\Gamma [u]  = u^{\prime 2}$.

\noindent
{\bf Example 3.} (Friedrich extension of a symmetric operator)

Let $D$ be a connected open set in $\RR^d$ with unit volume. Let $\PP=dx$ be the Lebesgue measure on $D$. Let $\Gamma$ be defined
on ${\cal C}_k^\infty(D)$ via
$$\Gamma[u,v]=\sum_{ij}\frac{\partial u}{\partial x_i}\frac{\partial v}{\partial x_j}a_{ij}, \quad\quad u,v\in{\cal C}_k^\infty(D)$$
where the functions $a_{ij}$ satisfy

$a_{ij}\in L^2_{loc}(D)\quad \frac{\partial a_{ij}}{\partial x_k}\in L^2_{loc}(D)\quad i,j,k =1,\ldots, d$,

$\sum_{ij}a_{ij}(x)\xi_i\xi_j\geq 0\quad\forall \xi\in D$,

$ a_{ij}(x)=a_{ji}(x) \quad \forall x\in D$,

\noindent then the pre-structure $(D,{\cal B}(D), \PP,{\cal C}_k^\infty(D), \Gamma)$ is closable.\\

Let us now come back to the question of using Dirichlet forms instead of second order vectors as germs of semi-martingales.\\

The first reason is the closedness property. That gives all the power to this theory. It is similar to $\sigma$-additivity in probability theory.
Without the closedness property, we have an apparently  more general framework (as additive 
set functions are more general than $\sigma$-additive ones), but  it becomes impossible to say anything on objects which are defined
by limits, error propagation is limited to explicit closed formulae. Instead, this closedness property allows to extend error calculus to
infinite dimensional frameworks and to propagate errors through typically limit objects as stochastic integrals. As David Hilbert argued against
intuitionists, more theorems is better. The philosopher Carl Popper made this mistake about axiomatization
of probability theory emphasing that his system (without $\sigma$-additivity) was more general than that of Kolmogorov (with $\sigma$-additivity).

What is particularly satisfying is that this closedness property is preserved by products. Any countable product of error structures is an error
structure and the theorem
on products  (cf. [Bou-03b]) gives explicitely the domain of the new $\Gamma$ operator. Starting 
with the Ornstein-Uhlenbeck structure in
dimension one, the infinite product of this structure by itself gives the Ornstein-Uhlenbeck structure on the Wiener space. 
Less surprisingly, the image of an error structure, defined in the most natural way, is still an error structure, as an image
of a probability space by a measurable map is still a probability space.\\

The second reason is related to simplicity. Let us come back to the first remark at the beginning of this part. We said  that 
starting with a centered
error, centeredness is lost after a non linear map. But what is preserved by image? Which property is an invariant ? It is the global property 
of symmetry with respect to a measure. If the operators describing the error are symmetric with respect to some measure, 
the image of the error has still this symmetry with respect to the image measure. Centeredness is nothing but symmetry with respect
to Lebesgue measure (not a probability measure, a $\sigma$-finite measure but this doesn't matter really here).\\

\noindent{\bf The gradient and the sharp ($\#$)}.

In addition to the operators $\Gamma$ and $A$ we will need the notion of {\it gradient} which is a linear (Hilbert valued) version of the 
standard deviation of the error.

\noindent{\bf Definition.} {\it Let ${\cal H}$ be a Hilbert space. A linear operator $D$ from $\DD$ into $L^2(\PP, {\cal H})$ is said to be a gradient (for $S$)
if $$\forall u\in\DD\quad \Gamma[u]=<Du,Du>_{\cal H}.$$}
A gradient always exists as soon the space $\DD$ is separable. It satisfies necessarily the chain rule :

\noindent{\bf Proposition} {\it Let $D$ be a gradient for $S$ with values in ${\cal H}$. Then $\forall u\in\DD^n, \forall F\in {\cal C}^1\cap Lip(\RR^n)$,
$$D[F\circ u]=\sum_{i=1}^n\frac{\partial F}{\partial x_i}\circ u D[u_i]\quad a.e.$$}

What we denote by the sharp $\#$ is a special case of gradient operator when ${\cal H}$ is chosen to be $L^2(\hat{\Omega},\hat{\cal A},\hat{\PP})$ where 
$(\hat{\Omega},\hat{\cal A},\hat{\PP})$ is a copy of $(\Omega ,\mathcal{A},\PP )$. It is particularly usefull for structures on the Wiener space
because stochastic calculus and Ito formula are available both on $(\Omega ,\mathcal{A},\PP )$ and $(\hat{\Omega},\hat{\cal A},\hat{\PP})$.

Let us give some definitions and notation we will need later on about the weighted Ornstein-Uhlenbeck structure on the 
Wiener space : let $B$ be a standard Brownian motion constructed as coordinates of the space ${\cal C}([0,1])$ equipped with the Wiener measure and let $\alpha$ be a
positive function in $L^1_{loc}[0,1]$, there exists an error structure
(cf. [Bou-03b]) satisfying
$$\Gamma[\int_0^1u(s)dB_s]=\int_0^1\alpha(s)u^2(s)ds$$
for $u\in{\cal C}([0,1])$. It is the mathematical expression of the following perturbation of the Brownian path :
$$\omega(s)=\int_0^s dB_u\mapsto \int_0^s e^{-\frac{\alpha(u)}{2}\varepsilon}dB_u+\int_0^s\sqrt{1-e^{-\alpha(u)\varepsilon}}d\hat{B}_u,$$
where $\hat{B}$ is an independent standard Brownian motion. This structure possesses the following $\#$-operator :
$$(\int_0^1u(s)dB_s)^\#=\int_0^1\sqrt{\alpha(s)}u(s)d\hat{B}_s,\quad \forall u\in L^2([0,1], (1+\alpha)dt),$$ which satisfies 
for regular adapted processes $H$
$$(\int_0^1 H_s dB_s)^\#=\int_0^1\sqrt{\alpha(s)}H_sd\hat{B}_s+\int_0^1 H_s^\# dB_s.$$

Let us end this part by a comment on  the passage from a random walk to the Brownian motion in the context
of erroneous quantities. Donsker's theorem says that if $U_n$ are i.i.d. square integrable centered random variables, the linear interpolation
of the random walk $\sum_{k=1}^n U_k$ i.e. the process
$$X_n(t)=\frac{1}{\sqrt{n}}\left(\sum_{k=1}^{[nt]}U_k+(nt-[nt])U_{[nt]+1}\right)$$
for $t\in[0,1]$, where $[x]$ denotes the entire part of $x$, converges in law on the space ${\cal C}([0,1])$ equipped by the 
uniform norm to a Brownian motion. Invariance principles follow giving a way to approximate properties of the Brownian motion 
by the corresponding ones of the 
random walk. A quite natural question is how this may be extended to the case where the $U_n$ are erroneous. To 
extend weak convergence of probability measures we use convergence of  Dirichlet forms  on
Lipschitz and ${\cal C}^1$ functions. Then supposing the errors on the $U_n$'s are equidistributed and uncorrelated, the error 
structure of the process $X_n$ converges to the Ornstein-Uhlenbeck structure on the Wiener space (cf. [Bou-05a]). Invariance 
principles follow giving approximations of the variance of the error of Brownian functionals, for example for the sup-norm 
of the paths :
$$\EE\Gamma[\|X_n(t)\|_\infty]=\EE\Gamma[\frac{1}{\sqrt{n}}\max_{1\leq k\leq n}|S_k|]\rightarrow
\EE[\int_0^1 (D_s[\|.\|_\infty])^2ds]=
\EE[{\cal T}]$$
where $D$ denotes the Ornstein-Uhlenbeck gradient with values in $L^2([0,1])$ and ${\cal T}$ is the random time 
where the absolute value of the Brownian 
path reaches its maximum.\\

\noindent {\sf\large III. Propagation of the error due to the Euler scheme.}\\

If an asset $X$ is represented by the solution of an sde, prices of options, hedging portfolios and other financial quantities 
are obtained by 
stochastic calculus as functionals of $X$. If we suppose the sde is solved using the Euler scheme, the asymptotic error on $X$ 
discovered by Kurtz and Protter in the spirit of a functional central limit theorem takes the form of a process 
solution to an other sde. In order to propagate this assymptotic error through stochastic calculus, we have to take the derivative in 
a suitable sense of non differentiable functionals as stochastic integrals. This may be performed by the theory of Dirichlet forms. 
Let us recall
the situation.\\

\noindent{\bf The error due to the Euler scheme.}

In 1991 Thomas Kurtz and Philipp Protter obtained an asymptotic estimate in law for the error
due to the Euler scheme 

In the simplest case, {considering the sde
$$X_t=x_0+\int_0^t a(X_s)dB_s+\int_0^tb(X_s)ds,$$
if $X_t^n$ is the Euler approximation of $X_t$ and $U^n=X^n-X$ then $(B,\sqrt{n}U^n)$ converges in law to $(B,U)$ where $U$ is solution to the linear sde

$$dU_t=a^{\prime}(X_t)U_tdB_t+b^{\prime}(X_t)U_tdt+\frac{1}{\sqrt{2}}a^{\prime}(X_t)a(X_t)dW_t,\quad\quad U_0=0,$$
where $W$ is a Brownian motion  independent of $B$.}

Such an ``extra-Brownian motion" appeared in a work of H. Rootzen [Roo-80] who studies limits of integrals of the form
$\int_0^t \psi_n(s)dB_s$ where $\psi_n$ is an adapted process. In the case where $\int_0^tf(B_s,s)dB_s$ is computed by the Euler scheme
$$\int_0^t\psi_n(s)dB_s=\sum_{i=0}^{[nt]}f(B_{\frac{i}{n}},i/n)(B_{\frac{i+1}{n}}-B_{\frac{i}{n}})+f(B_{\frac{[nt]}{n}},[nt]/n)(B_t-B_{\frac{[nt]}{n}})$$
he obtains for regular $f$ 
$$\sqrt{n}\left(\int_0^. \psi_ndB-\int_0^.f(B_s,s)dB_s\right)\stackrel{d}{\Rightarrow}\frac{1}{\sqrt{2}}\int_0^.f^{\prime}_x(B_s,s)dW_s.$$
This kind of result is restricted to adapted approximations. 
As  Wong and Zakai have shown (1965 ) other natural  approximations of the brownian motion give rise to stochastic integrals in the sense of Stratonowitch

The discovery of the asymptotic error due to the Euler scheme has been followed by a series of works
which extend it to the case of an sde with respect to a continuous or discontinuous semi-martingale
and which obtain some statements as necessary and sufficient conditions ([Ja-Pr-98], [Ja-Ja-M\'e-03]).

In addition, asymptotic expansions have been recently obtained
by the stochastic calculus of variation [Ma-Th-03].

In the sequel, we shall consider the result of   Kurtz-Protter in dimension 1 under  the following  form:

{Let $X_t$ be the  solution starting at  $x_0$ to the sde
$$dX_t=a(X_t,t)dB_t+b(X_t,t)dt,$$
let $X_t^n$ be the approximate solution obtained by the Euler method, which may be written 
$$X_0^n=x_0;\quad\quad\quad dX_t^n=a(X^n_{\frac{[nt]}{n}},[nt]/n)dB_t+b(X^n_{\frac{[nt]}{n}},[nt]/n)dt$$
and let $U_t^n=X^n_t-X_t$ be the approximation error, then if $a$ and $b$ are ${\cal C}^1$ with linear growth
$$(B,\sqrt{n} U^n)\stackrel{d}{\Rightarrow} (B,U)\hspace{2cm}{\mbox{ on }}\;\;{\cal C}([0,1])$$
where the process $U$ may be represented as 
$$U_0=0\quad\quad\quad dU_t=a^{\prime}_x(X_t,t)U_tdB_t+b^{\prime}_x(X_t,t)U_tdt+\frac{1}{\sqrt{2}}a^{\prime}_x(X_t,t)a(X_t,t)dW_t$$}
which is solved by the usual method of variation of the constant : introducing the process
$$M_t=\exp\left\{\int_0^ta^{\prime}_x(X_s,s)dB_s-\frac{1}{2}\int_0^ta^{\prime 2}_x(X_s,s)ds+\int_0^tb^{\prime}_x(X_s,s)ds\right\}$$
gives
$$U_t=M_t\int_0^t\frac{a(X_s,s)a^{\prime}_x(X_s,s)}{\sqrt{2}M_s}dW_s.$$

Let us consider the weighted Ornstein-Uhlenbeck error structure on the Wiener space with weight $\alpha$ as explain above. If the  coefficients $a$ and $b$ are regular, then  $X_t\in\DD$ and $X_t^{\#}$ satisfies\\

$\displaystyle(*)\quad X_t^{\#}=\int_0^t \!a^\prime_x(X_s,s)X_s^{\#}dB_s+\int_0^t \!a(X_s,s)\sqrt{\alpha(s)} d\widehat{B_s}+\int_0^t\!b^\prime_x(X_s,s)X_s^{\#}ds
$\\

Comparing with the equation of the asymptotic error due to the Euler scheme\\

$\displaystyle(**)\quad U_t=\int_0^t \!a^\prime_x(X_s,s)U_sdB_s+\int_0^t \!a(X_s,s)\frac{a^\prime_x(X_s,s)}{\sqrt{2}}
 dW_s+\int_0^t\!b^\prime_x(X_s,s)U_sds
$\\

\noindent shows that

- if we could take a random and adapted weight 
$\alpha(t)=\frac{1}{2}a^{\prime 2}_x(X_t,t)$

- if the obtained structure is  closable 
with carr\'e du champ and if the calculus of the  $\#$-operator is still (*)

\noindent then  $X^{\#}$ would be the asymptotic error due to the Euler scheme, and we would be able to propagate this error through
the stochastic computations   obtaining the variance of the error on any r. v. $Y\in\DD$ by the equation
 $\Gamma[Y]=\widehat{\EE}[Y^{\# 2}]$. \\

\noindent{\bf The Ornstein-Uhlenbeck structure with random weight.}

From now on $\alpha$ is a measurable random process defined on the Wiener space, non negative, non necessarily adapted. We assume
that this process satisfies 
$\EE\int_0^1\alpha_t dt<+\infty$, and 
$\alpha(\omega, t)\geq k(t)>0\quad \PP\times dt{\mbox{-a.e.}}$ where $k$ is deterministic.

{Let us denote $\DD_{ou}^k$ the domain of the Ornstein-Uhlenbeck structure with deterministic weight $k$ and  $D_{ou}^k$
  its gradient.
 On the domain
$$\DD=\left\{ Y\in\DD_{ou}^k :\quad \int_0^1\EE[(D_{ou}^k[Y](t))^2\frac{\alpha(t)}{k(t)}] dt<+\infty\right\}$$
which is dense, the form
$${\cal E}[Y]=\frac{1}{2}\int_0^1\EE[(D_{ou}^k[Y](t))^2\frac{\alpha(t)}{k(t)}] dt$$
is Dirichlet and admits
$$\Gamma[Y]=\int_0^1(D_{ou}^k[Y](t))^2\frac{\alpha(t)}{k(t)} dt$$ as carr\'e du champ operator.}

Indeed,  let ${\cal V}$ be the space of linear combinations of exponentials of the form
$Y=\exp\{i\int_0^1 h_u dB_u\}$ with $h$ deterministic bounded, 
by  $\int_0^1\EE\alpha(t) dt<+\infty$, we have ${\cal V}\subset\DD$ and 
$D_{ou}^k[Y]=Y(ih\sqrt{k})$ hence $\DD$ is dense.

Let $X_n$ be a Cauchy sequence in $L^2$ and for ${\cal E}$. Let $X$ be the limit of $X_n$ in $L^2$. 
Then $X_n$ is Cauchy for ${\cal E}_{ou}^k$ which is closed, hence $X\in\DD_{ou}^k$ and there exists a sub-sequence
 $X_{n^\prime}$ such that
$$D_{ou}^k[X_{n^\prime}]\rightarrow D_{ou}^k[X]\quad\quad \EE\times dt{\mbox{-p.s.}}$$
and by Fatou's lemma
$$\int_0^1\EE[(D_{ou}^k[X])^2\frac{\alpha(t)}{k(t)}]dt=$$
$$=\int_0^1\EE[\lim(D_{ou}^k[X_{n^\prime}])^2\frac{\alpha(t)}{k(t)}]dt\leq 
\liminf\int_0^1\EE[(D_{ou}^k[X_{n^\prime}])^2\frac{\alpha(t)}{k(t)}]dt<+\infty$$ since $X_n$ is Cauchy for ${\cal E}$. Hence $X\in\DD$. Now again by the
Fatou's lemma we show as classically that
 $X_n$ converges to $X$ in $\DD$.

Contractions operate on $({\cal E},\DD)$ by the functional calculus for $D_{ou}^k$ hence $({\cal E},\DD)$ is a Dirichlet form. The definition 
of the  
carr\'e du champ operator (def 4.1.2 of [Bo-Hi-91]) is satisfied.

 {The generator $(A,{\cal D}A)$ is given by
$$\begin{array}{rl}
{\cal D}A=&\{F\in\DD\quad\exists G\in L^2\quad\forall H\in\DD\quad
\frac{1}{2}\EE\int_0^1D_{ou}^k[F]D_{ou}^k[H]\frac{\alpha(t)}{k(t)}dt
\left.=-<G,H>\right\}\\
AF=&G
\end{array}$$
\noindent{hence if $F\in{\cal D}A$ then $\frac{\alpha(t)}{k(t)}D_{ou}^k[F]\in{\mbox{dom}}\delta_{ou}^k$ and
$$AF=-\frac{1}{2}\delta_{ou}^k[\frac{\alpha}{k}D_{ou}^kF].$$}
where $\delta_{ou}^k$ is the  Skorokhod integral with weight $k$.\\

\noindent{\bf Adapted case.}

Let us now add the hypothesis that $\alpha$ is adapted. 
 { If $h$ is in $L^\infty(\RR_+)$
$$\EE\Gamma[F,\int_0^1 hdB]]=\EE[F\int_0^1 h(s)\alpha(s)dB_s].$$}
 {If $F,G\in\DD\cap L^\infty$
$$\EE[G<DF,h\sqrt{\alpha}>]=-\EE[F<DG,h\sqrt{\alpha>}]+\EE[FG\int h\alpha dB].$$}
 {And if $v$ is  adapted and in $\mbox{dom}\delta$
$$\delta[v]=\int_0^1v_s\sqrt{\alpha_s}dB_s.$$}
 At last, for finance, the following properties are important, they use  the fact that $\alpha$ is adapted
 {
$$A[\EE[X|{\cal F}_s]]=\EE[A^s[X]|{\cal F}_s]$$ where $A^s$ is constructed as $A$ with the weight $\alpha(t)1_{\{t\leq s\}}$,
$$D[\EE[X|{\cal F}_s]](t)=\EE[D[X](t)1_{t\leq s}|{\cal F}_s]$$
$\EE[.|{\cal F}_s]$ is an orthogonal projector in $\DD$
$$(\EE[X|{\cal F}_s])^{\#}=\EE[X^{\#_s}|{\cal F}_s]$$
where $\#_s$ is constructed as $\#$ with the weight $\alpha(t)1_{\{t\leq s\}}$.
If $X$ is ${\cal F}_t$-measurable, then  $AX$, $\Gamma[X]$ are ${\cal F}_t$-measurables.
}

 {Concerning the operator $\#$
we have the formulae
$$\left(\int_0^1 \xi_s\,dB_s\right)^{\#}=\int_0^1\xi_s^{\#}\,dB_s+\int_0^1\xi_s\sqrt{\alpha_s}d\widehat{B}_s$$}
Hence formula ($\star$) is satisfied.\\

\noindent{\bf Application to diffusion models.}

Let us consider the following model of an asset 
$$dX_t=X_t\sigma(X_t,t)dB_t+X_t r(t)dt$$
and let us put on the Wiener space the Ornstein-Uhlenbeck structure with weight 
$$\alpha_t=\frac{a^{\prime 2}(X_t,t)}{2}=\frac{(\sigma(X_t,t)+X_t\sigma^\prime_x(X_t,t))^2}{2}$$
which represents the asymptotic error due to the Euler scheme. $\sigma$ is supposed to be strictly positive, ${\cal C}^1$and Lipschitz 
and the preceding hypotheses on $\alpha$ are assumed.

Such a modelisation is coherent.  The error is attached to the asset $X$
and any functional of $X$, including the Brownian motion itself and its error may be computed thanks to the equation
$$dB_t=\frac{dX_t}{X_t\sigma(X_t,t)}-X_t r(t) dt$$
which  gives
$$(B_t)^\#=\int_0^t \sqrt{\alpha(s)}d\widehat{B}_s\quad\quad\quad \Gamma[B_t]=\int_0^t \alpha(s)ds.$$
Let us show how  financial calculi may be performed before proposing some comments on the use of such an analysis. Puting
$M_t=\exp\{\int_0^t\sqrt{\alpha_s}dB_s-\frac{1}{2}\int_0^t\alpha_sds+\int_0^tr(s)ds\}$ we have
$$\Gamma[X_t]=M_t^2\int_0^t\frac{X_s^2\sigma^2(X_s,s)}{M_s^2}\alpha_sds$$
$$\Gamma[X_s,X_t]=M_sM_t\int_0^{s\wedge t}\frac{X_u^2\sigma^2(X_u,u)}{M_u^2}\alpha_udu.$$
The price of a European option with payoff $f(X_T)$ at exercise time $T$
$$V_t=\EE[(\exp-\!\int_t^T\!\!r(s)ds)f(X_T)|{\cal F}_t]$$
becomes erroneous (in the sense of error structures) with an error obtained thanks to the $\#$:
$$\Gamma[V_t]=(\exp-2\!\!\int_t^T\!\!\!r(s)ds)(\EE[f^\prime(X_T)M_T|{\cal F}_t])^2\frac{\Gamma[X_t]}{M_t^2}$$
$$
\Gamma[V_s,V_t]=\left(\exp(-\!\!\int_s^T\!\!\!r(u)du-\!\!\int_t^T\!\!\!r(v)dv)\right)
\EE[f^\prime(X_T)M_T|{\cal F}_s]\EE[f^\prime(X_T)M_T|{\cal F}_t]\frac{\Gamma[X_s,X_t]}{M_sM_t}
$$
 The quantity of asset in the hedging portfolio is 
$$H_t=(\exp-\!\int_t^T\!\!r(s)ds)\EE[f^\prime(X_T)M_T|{\cal F}_t]\frac{1}{M_t}$$
and we have
$$\Gamma[H_t]=(\exp-2\!\!\int_t^T\!\!\!r(s)ds)(\EE[\frac{M_T}{M_t}(f^{\prime\prime}(X_T)M_T+f^\prime(X_T)Z_t^T)|{\cal F}_t])^2\frac{\Gamma[X_t]}{M_t^2}$$
with
$$Z_t^T=\int_t^TL_sdB_s-\int_t^T\sqrt{\alpha_s}L_sM_sds$$
$$L_s=a^{\prime\prime}_{x^2}(X_s,s)=2\sigma^\prime_x(X_s,s)+X_s\sigma^{\prime\prime}_{x^2}(X_s,s)$$
 It is still true, as in the case of deterministic weight (cf [Bou-03b]),  that the proportional error on  $X_t$ divided by the volatility :
$$\frac{\sqrt{\Gamma[X_t]}}{X_t}\cdot\frac{1}{\sigma(X_t,t)}$$ is a finite variation process 
(cf [Ba-Ma-Ma-Re-Th-03]) on the``feed back" effect).\\

\noindent{\bf Discussion.}

Thanks to this construction of an error structure, i.e. a local Dirichlet form with squared field operator, on the Wiener space, hence
by image, on ${\cal C}([0,1])$ equipped with the law of the process $X$, we have at our disposal  a powerful mean 
  to propagate
the error done on $X$ toward sufficiently smooth functionals of $X$. In order to assess the interest of this tool, the 
question arises of knowing whether the propagated error is the same as the one we would obtain by a direct computation
of the functional thanks to the approximation $X^n$ of $X$. For instance, in the simplest case, does the convergence in law
$$\sqrt{n}(f(X_t^n)-f(X_t))\stackrel{d}{\Rightarrow} f^\prime(X_t)X_t^{\#}$$
hold for $f\in{\cal C}^1\cap Lip$ ? Can we justify  {\it an asymptotic calculus principle} which says that the Dirichlet form allows
effectively to compute the errors on the quantities which are erroneous because of the approximation $X^n$ of $X$ ? We will not
exhaustively examine this principle here, for it is a too large enterprise. Nevertheless, in the important current of research 
whose fruitfulness
has been confirmed these last twenty years, which may be called the ``tightness programm", the authors, among which we must at least
quote P.-A. Meyer, W. A. Zheng, J. Jacod, A. N. Shiryaev, A. Jakubowski, J. M\'emin, G. Pag\`es, T. G. Kurtz, P. Protter, 
L. S\l omi\'nski, D. Talay, V. Bally, A. Kohatsu-Higa and many others, have already done a major part of the work by stating their 
results of convergence in law, of stable convergence, of tightness of processes, under a sufficiently general form 
for propagating iteratively the properties through stochastic integrals and sde's in the semi-martingale framework.

Let us give some results in the direction of this {\it asymptotic calculus principle} keeping the hypotheses of the present part III.

Let $F$ be a real function of class ${\cal C}^1$ and Lipschitz defined on ${\cal C}([0,1])$ equipped with the uniform norm. Such a function
satisfies
$$F(x+h)=F(x)+<F^\prime(x),h>+\|h\|\varepsilon_x(h)\quad\quad\forall x,h\in{\cal C}([0,1])$$
where the mapping $x\mapsto F^\prime(x)$ is continuous and bounded with values in the Banach space of Radon 
measures on $[0,1]$, $\varepsilon_x(h)$ is bounded in $x$ and $h$, and goes to zero when $h\rightarrow 0$ in ${\cal C}([0,1])$.
Then we have
$$\sqrt{n}(F(X^n)-F(X))\stackrel{d}{\Rightarrow}  (F(X))^\#=\int_{[0,1]}X_t^\#\;F^\prime(X)(dt).$$
{\it Proof.} The equality in the right hand side comes from the functional calculus in error structures (see [Bou-05a]). Puting
$U^n=X^n-X$ as before, the fact that $\sqrt{n}\|U^n\|\varepsilon_X(U^n)$ tends to zero in probability, reduces the proof to the study
of the convergence in law of 
$$<F^\prime(X),\sqrt{n}U^n>=\sqrt{n}\int(X^n_t-X_t)\;F^\prime(X)(dt)$$
to $\int X_t^\#\;F^\prime(X)(dt)$. Considering the measure $F^\prime(X)(dt)$ as the differential of a finite variation process adapted
to the constant filtration ${\cal G}_t={\cal B(C}(]0,1]))$, the fact that the process to be integrated $\sqrt{n}U^n$ 
converges stably to $X^\#$ implies (cf [Ku-Pr-91b]
thm 2.2)
that the stochastic integral $\int \sqrt{n}U^n\;F^\prime(X)(dt)$ converges in law to $\int X_t^\#\;F^\prime(X)(dt)$.

We obtain also the convergence in law of the stochastic integrals $ H.\sqrt{n}U^n\stackrel{d}{\Rightarrow}  H.X^\#$ for $H$ deterministic or adapted
and that of 
$$\sqrt{n}(\int_0^1 f(X^n_s,s)dX_s^n-\int_0^1f(X_s,s)dX_s)$$
to 
$$(\int_0^1f(X_s,s)dX_s)^\#=\int_0^1f^\prime(X_s,s)X^\#_sdX_s+\int_0^1f(X_s,s)dX_s^\#$$
for $f$ ${\cal C}^1$ and Lipschitz.\\

More generally, we can make more explicit the research programm of determining {\it the domain} of the asymptotic  calculus.

Let $X_n$ and $X$ be two random variables with values in a measurable set $(E,{\cal F})$, and let $\alpha_n$ be a 
sequence of positive numbers. Let ${\cal D}_0$ denote a set of {\it simple} functions included in $L^2(\PP_X)$ and in 
$L^2(\PP_{X_n})$ $\forall n$. Let us suppose that there exists an error structure 
$$S=(E,{\cal F},\PP_X,\DD,\Gamma)$$
such that ${\cal D}_0\subset \DD$ and $\forall\varphi\in{\cal D}_0$
\begin{equation}\lim_n \alpha_n\EE[(\varphi(X_n)-\varphi(X))^2]=\EE[\Gamma[\varphi]]\end{equation}
we shall say that the asymptotic calculus principle extends to ${\cal D}$ for ${\cal D}_0\subset{\cal D}\subset\DD$ 
if the limit (1) extends to $\psi\in{\cal D}$.

If, as above, a $\#$-operator is available (which occurs as soon as $\DD$ is separable), in order to prove (1) on ${\cal D}$, since 
$\#$ is a closed operator, it suffices for any $\psi\in{\cal D}$ to find a sequence $\varphi_p\in{\cal D}_0$ such that
$$
\begin{array}{l}
i)\quad \varphi_p\rightarrow \psi \mbox{ in }L^2(\PP_X)\\
ii)\quad \varphi_p^\# {\mbox{ converges in }} L^2(\PP_X\times\widehat{\PP_X})\\
iii)\quad \alpha_n\EE[\psi(X_n)-\psi(X))^2]{\mbox{ may be approximated uniformly in {\it n} by }}\\
\qquad\alpha_n \EE[\varphi_p(X_n)-\varphi_p(X))^2].
\end{array}
$$
When $(E,{\cal F})$ is a normed vectorspace, obtaining (1) from a convergence in law of 
$\sqrt{\alpha_n}(\varphi(X_n)-\varphi(X))$ uses generally
a uniform integrability of $\alpha_n\|X_n-X\|^2$. We shall go deaper in this problem in a separate work. \\

Let us end by some remarks from the point of view of finance. The interest of considering a financial asset as erroneous is not evident
since it is one of the best known quantities continuously quoted in a financial market. Such an error may be justified (cf [Bou-03b])
by the inaccuracy of the instants of transaction, possibly also  to represent an infinitesimal bid-ask. But this would rather 
justify specifically constructed error structures instead of the one induced by the Euler scheme. This error structure is relevant only in order
to assess the errors in Monte Carlo simulations performed to calculate financial quantities in a given model.

Several authors ([Du-Pr-89], [Ha-My-05]) remarked that the stochastic integral which is the active hedge
of a future contingent claim, in a model where the underlying asset is a semi-martingale, is an instance of application of limit theorems on
discretization errors. This is different from the Euler scheme error and it would be worth to examine this error from the point 
of view of an asymptotic Dirichlet form.

A more general and complete study of the bias operators and the Dirichlet form yielded by an approximation, 
with applications related to the part III of the present paper, is to appear, ([Bou-06]).

\begin{list}{}
{\setlength{\itemsep}{0cm}\setlength{\leftmargin}{0.5cm}\setlength{\parsep}{0cm}\setlength{\listparindent}{-0.5cm}}
  \item\begin{center}
{\small REFERENCES}
\end{center}\vspace{0.4cm}

[Al-Ma-Nu-00] {\sc Alos E., Mazet O., Nualart D.} ``Stochastic calculus with respect to fractional Brownian motion with Hurst parameter lesser
 than 1/2" {\it Stochastic Process. Appl.} 86, (2000).

[Ba-Ta-96a] {\sc Bally V., Talay D.} ``The law of the Euler scheme for stochastic differential equations : I. Convergence rate of the distribution function", 
{\it Prob. Th. and Rel. Fields} vol 2 No2, 93-128 (1996)

[Ba-Ta-96b] {\sc Bally V., Talay D.} ``The law of the Euler scheme for stochastic differential equations : II. Convergence rate of the density", 
{\it Monte Carlo Methods and Appl.} vol 104, No1, 43-80 (1996)

[Ba-Ma-Ma-Re-Th-03] {\sc Barucci E., Malliavin P., Mancino M.E., Ren\'o R., Thalmaier A.} ``The price volatility feedback rate: an implementable indicator of market stability" 
{\it Math. Finance}, 13, 17-35, 2003.

[Bo-Ek-To-04] {\sc Bouchard B., Ekeland I., Touzi N.}``On the Malliavin approach to Monte Carlo approximation 
of conditional expectations", {\it Finance Stoch.} 8, 45-71, (2004),

[Bou-85] {\sc Bouleau N.} ``About stochastic integrals with respect to processes which are not semi-martingales" {\it Osaka J. Math.}
22, 31-34, (1985),

[Bou-03a] {\sc Bouleau N.} ``Error calculus and path sensitivity in Financial models", 
{\it Mathematical Finance} vol 13/1, 115-134, (2003).

[Bou-03b] {\sc Bouleau N.} {\it Error Calculus for Finance and Physics, the Language of Dirichlet Forms}, De Gruyter, 2003.

[Bou-04] {\sc Bouleau N.} {\it Financial Markets and Martingales, Observations on Science and Speculation}, Springer 2004.

[Bou-05a] {\sc Bouleau N.} ``Th\'eor\`eme de Donsker et formes de Dirichlet" {\it Bull. Sci. Math. }129, (2005), 369-380.

[Bou-05b] {\sc Bouleau N.} ``Improving Monte Carlo simulations by Dirichlet forms" C. R. Acad. Sci. Paris Ser I (2005)

[Bou-06] {\sc Bouleau N.} ``When and How an error yields a Dirichlet form" {\it J. of Functional Analysis}, to appear.

[Bo-Hi-91] {\sc Bouleau N., Hirsch F.}  {\it Dirichlet Forms and Analysis on Wiener Space,} De Gruyter, (1991).

[Ca-Fe-Nu-98] {\sc Caballero M. E., Fernandez B., Nualard, D.} ``Estimation of densities and applications" {\it J. of Theoretical Prob.} 11,
No3, (1998),

[De-Gr-Po-99] {\sc Denis l., Grorud A., Pontier M.} ``Formes de Dirichlet sur un espace de Wiener-Poisson, application au grossissement de 
filtration" {\it S\'em. Prob. XXXIV} Lect. N. in Math. 1729, Springer (1999),

[Du-Pr-89] {\sc Duffie D., Protter P.,}``From discrete to continuous time finance : weak convergence of the financial gain process" (1989)
(unpublished, cited in [Ku-Pr-91b])

[El-Pr-04] {\sc El Kahtib Y., Privault N.} ``Computation of Greeks in a market with jumps via Malliavin calculus" {\it Finance and Stoch.} 8, 161-179, (2004),

[Eme-89] {\sc Emery M.} {\it Stochastic calculus on Manifolds} Springer (1989),

[Er-Ru-98] {\sc Errami M., Russo F.} ``Covariation de convolution de martingales" {\it C. R. Acad. Sci. Paris}, s1, 326, 601-609, (1998),

[Fe-LaP-99] {\sc Feyel D., La Pradelle A. de},``On fractional Brownian processes" {\it Potential Anal.} 10(3), 273-288, (1999),

[Fl-Ru-Wo-00] {\sc Flandoli F., Russo F., Wolf J.} ``Some SDEs with distributional drift, I. General calculus" {\it Osaka J. Math.}
40, No2, 493-542, (2003).

[F\"{o}l-80] {\sc F\"ollmer H.} ``Dirichlet processes" {\it in Stochastic Integrals, lect. Notes in Math.} No 851, 476-478, (1981).

[F\"{o}l-81] {\sc F\"ollmer H.} ``Calcul d'Ito sans probabilit\'e" {\it in S\'em. Prob. XV}, Lect. N. in Math. 850,  143-150,  Springer (1981),

[F\"o-Pr-Sh-95] {\sc F\"ollmer H., Protter P., Shiryaev A. N.} ``Quadratic covariation and an extension of Ito's formula" {Bernoulli} 1, 149-169,
(1995),

[Fo-La-Le-Li-To-99] {\sc Fourni\'e E., Lasry J. M., Lebuchoux J., Lions P. L., Touzi N.} ``Applications of Malliavin calculus to Monte Carlo
methods in finance" {\it Finance and Stoch.} 3, 391-412, (1999),

[Fo-La-Le-Li-01] {\sc Fourni\'e E., Lasry J. M., Lebuchoux J., Lions P. L.}``Applications of Malliavin calculus to Monte Carlo
methods in finance II" {\it Finance and Stoch.}(2001),

[Fuk-80] {\sc Fukushima M.} {\it Dirichlet forms and Markov processes}, North-Holland (1980),

[Fu-Os-Ta-94] {\sc Fukushima M., Oshima Y., Takeda M.} {\it Dirichlet forms and symmetric Markov processes}, De Gruyter 1994.

[Go-Ko-01] {\sc Gobet E., Kohatsu-Higa A.} ``Computation of Greeks for barrier and lookback options using Malliavin calculus" {\it Electron. Com. in Prob.} 8, 51-62, (2003).

[Gr-No-03] {\sc Gradinaru M., Nourdin I.} ``Approximation at first and second order of  the $m$-variation of the fractional Brownian motion"
{\it Electron. Com. Prob.} 8, 1-26, (2003)

[Gr-Ru-Va-03] {\sc Gradinaru M., Russo F., Vallois P.} ``Generalized covariation, local time and Stratonowich-Ito's formula for fractional Brownian
motion with Hurst index$\geq$ 1/4" {\it Annals of Prob.} 31, No4, 1772-1820, (2003).

[Ha-My-05] {\sc Hayashi T., Mykland P.A.} ``Evaluating hedging errors : an asymptotic approach" 
{\it Math. Finance}, vol 15, No 2, 309-343, (2005)

[Jac-85] {\sc Jacod J.,} ``Th\'eor\`emes limites pour les processus" {\it Lect. Notes Math. } vol 1117, Springer 1985.

[Ja-Ja-M\'e-03] {\sc Jacod J., Jakubowski A., M\'emin J.} ``About asymptotic errors in discretization processes" {\it Ann; of Prob.}31, 592-608, (2003),

[Ja-Pr-98] {\sc Jacod J., Protter P.} ``Asymptotic error distributions for the Euler method for stochastic 
differential equations'' {\it Ann. Probab.} 26, 267-307, (1998)

[Ja-Sh-87] {\sc Jacod J., Shiryaev A.N.,} {\it Limit Theorems for Stochastic Processes,} Springer, 1987.

[Ja-M\'e-Pa-89] {\sc Jakubowski A., M\'emin J., Pag\`es G.} ``Convergence en loi des suites d'int\'egrales stochastiques sur l'espace de Skorokhod" {\it 
Probab. Th. Rel. Fields} 81, 111-137, 1989.

[Ko-Pe-02] {\sc Kohatsu-Higa A., Pettersson R.} ``Variance reduction methods for simulation of densities on Wiener space", 
{\it SIAM J. Numer. Anal.} Vol 40, No2, 431-450, (2002)

[Ku-Pr-91a] {\sc Kurtz T.; Protter P.} ``Wong-Zakai corrections, random evolutions and simulation schemes for SDEs" {\it Stochastic Analysis} 331-346,
Acad. Press, 1991.

[Ku-Pr-91b] {\sc Kurtz, T.; Protter, Ph.} ``Weak limit theorems for stochastic integrals ans stochastic differential equations" {\it Ann. Probab.} 19, 1035-1070, 1991.

[Le-Na-Nu-03] {\sc Leon J. A., Navarro R., Nualart D.} ``An anticipating calculus approach to the utility maximization
of an insider" {\it Math. Finance} 13, No1, 171-185, (2003).

[LeJ-78] {\sc Le Jan Y.} ``Mesures associ\'ees \`a une forme de Dirichlet, applications" {\it Bull. Soc. Math. France} 106, 61-112, (1978),

[Ly-Zh-98] {\sc Lyons T., Zheng W.} ``A crossing estimate for the canonical process on a Dirichlet space and tightness result" 
{\it in Colloque Paul L\'evy, Ast\'erisque} No 157-158, pp 249-271, (1998),

[Ma-R\"o-92] {\sc Ma Z. M., R\"ockner M.} {\it Introduction to the Theory of (Non-Symmetric) Dirichlet Forms} Springer (1992),

[Ma-Th-03] {\sc Malliavin P., Thalmaier A.} ``Numerical error for SDE: Asymptotic expansion and hyperdistributions", {\it C. R. Acad. Sci. Paris}
ser. I 336 (2003) 851-856

[Ma-Th-05] {\sc Malliavin P., Thalmaier A.} {\it Stochastic Calculus of Variations in Mathematical Finance}, Springer, (to appear 2005),

[Me-Pr-03] {\sc Mensi M., Privault N.} ``Conditional calculus and enlargement of filtration on Poisson space" {\it Stoch. Anal. and Appl.} 21, 183-204,(2003),

[Mey-80] {\sc Meyer P.-A.}``G\'eom\'etrie diff\'erentielle stochastique" {\it in S\'em. Prob. XVI suppl. } 165-207, Lect. 
N. in Math. 921, Springer (1982),

[Nua-95] {\sc Nualart N.}  {\it The Malliavin calculus and related topics}. Springer, 1995.

[Osh-92] {\sc Oshima, Y} ``On a construction of Markov processes associated with time dependent Dirichlet spaces"
 {\it Forum Math.} No 4, 395-415, (1992).

[Roo-80] {\sc Rootz\'en, H.} ``Limit distribution for the error in approximation of stochastic integrals''
 {\it Ann. Probab. } 8, 241-251, (1980).

[Ru-Va-95] {\sc Russo F., Vallois P.} ``The generalized covariation process and It\^o formula" {\it Stochastic Proc. Appl.} 59, 81-104, (1995),

[Ru-Va-96] {\sc Russo F., Vallois P.}  ``It\^o formula for $C^1$-functions of semi-martingales" {\it Prob. Th. Rel. Fields} 104, 27-41, (1996),

[Ru-Va-00] {\sc Russo F., Vallois P.} ``Stochastic calculus with respect to  a finite quadratic variation process" {\it Stochastics and Stoch. Rep.} 70, 1-40, (2000),

[Ru-Va-Wo-01] {\sc Russo F., Vallois P., Wolf J.} ``A generalized class of Lyons-Zheng processes" {\it Bernoulli} 7, No2, 363-379, (2001),

[Sch-82] {\sc Schwartz L.}``G\'eom\'etrie diff\'erentielle du 2\`eme ordre, semi-martingales et \'equations diff\'erentielles stochastiques 
sur une vari\'et\'e diff\'erentielle" {\it in S\'em. Prob. XVI suppl. } 1-150, Lect. N. in Math. 921, Springer (1982),

[S\l o-89] {\sc S\l omi\'nski, L.} ``Stability of strong solutions of stochastic differential equations'' 
{Stochastic Process. Appl.} 31, 173-202, (1989).

[Sta-99] {\sc Stannat, W.} ``The theory of generalized Dirichlet forms and its applications in analysis and stochastics" {\it Mem. 
Amer. Math. Soc.} 142, No 678, (1999).

[Tru-00] {\sc Trutnau G.} ``Stochastic calculus of generalized Dirichlet forms and applications to stochastic differential equations in infinite
dimensions" {\it Osaka J. Math.} 37(2), 315-343, (2000),

[Wol-97] {\sc Wolf J.} ``An Ito formula for Dirichlet processes" {\it Stochastics and Stoch. Rep.} 62(2), 103-115, (1997),

[Wo-Za-65] {\sc Wong E., Zakai M.} ``On the convergence of ordinary integrals to stochastic integrals" 
{\it Ann. Math. Statist.} 36, 1560-1564, (1965)

[Zah-98] {\sc Zahle M.} ``Integration with respect to fractal functions and stochastic calculus" 
{\it Prob. Th. Rel. Fields} 21, 333-374, (1998),

\end{list}

\end{document}